\documentclass[a4paper,11pt,reqno]{amsart}
\usepackage[centering,includeheadfoot,margin=3cm]{geometry}
\usepackage{graphics,enumitem,epsfig,textcomp}
\usepackage{amsfonts,amsmath,amssymb,euscript,color,mathrsfs}
\usepackage[utf8]{inputenc}	

\usepackage[caption=false,subrefformat=parens,labelformat=parens]{subfig}

\usepackage{url}

\usepackage[symbol]{footmisc}

\usepackage{hyperref}

\usepackage{lineno}

\numberwithin{equation}{section}

\newtheorem{theorem}{Theorem}[section]

\newtheorem{lemma}{Lemma}[section]

\def\d{\,\mathrm{d}}

\def\N{\mathbb{N}}
\def\R{\mathbb{R}}
\def\C{\hbox{\rlap{\kern.24em\raise.1ex\hbox
      {\vrule height1.3ex width.9pt}}C}}
\def\P{\hbox{\rlap{I}\kern.16em P}}
\def\Q{\hbox{\rlap{\kern.24em\raise.1ex\hbox
      {\vrule height1.3ex width.9pt}}Q}}
\def\M{\hbox{\rlap{I}\kern.16em\rlap{I}M}}
\def\Z{\hbox{\rlap{Z}\kern.20em Z}}
\def\({\begin{eqnarray}}
\def\){\end{eqnarray}}
\def\[{\begin{eqnarray*}}
\def\]{\end{eqnarray*}}
\def\part#1#2{\frac{\partial #1}{\partial #2}}

\def\Norm#1{\left\| #1 \right\|}
\def\pmb#1{\setbox0=\hbox{$#1$}
  \kern-.025em\copy0\kern-\wd0
  \kern-.05em\copy0\kern-\wd0
  \kern-.025em\raise.0433em\box0 }
\def\bar{\overline}

\def\d{\,\mathrm{d}}
\def\N{\mathbb{N}}
\def\R{\mathbb{R}}

\def\epsilon{\varepsilon}

\def\P{\mathbb{P}}
\def\QQ{\mathcal{Q}}
\def\ZZ{\mathcal{Z}}

\def\I{\mathcal{I}}
\def\intI{\int_{\I}}
\def\iintII{\intI\intI}
\def\II{\square}

\def\dd{\bar \nabla}
\def\II{\I\times\I}
\def\L{\mathrm{L}}

\newcommand\GraphG{\mathbb{G}}
\newcommand\Vset{\mathbb{V}}
\newcommand\Eset{\mathbb{E}}

\def\C{\mathfrak{c}}

\def\Fcont{\mathcal{F}}

\def\b{b}  
\def\barb{\tilde\b}  

\def\Fkin{\Fcont_{\mathrm{kin}}}

\def\Fmet{\Fcont_{\mathrm{met}}}

\def\wto{\rightharpoonup}

\def\Ltwozero{L^2_0}
\def\LtwozeroN{L^2_{0,N}}
\def\BB{\mathbb{B}}

\begin{document}


\centerline{\Large\textbf{Rigorous dense graph limit of a model for}}
\centerline{\Large\textbf{biological transportation networks}}
\vskip 7mm

\centerline{
	{\large Nuno J. Alves}\footnote{Faculty of Mathematics, University of Vienna,
	    Oskar-Morgenstern-Platz 1, 1090 Wien, Austria;
		{\it nuno.januario.alves@univie.ac.at}}\qquad and\qquad
	{\large Jan Haskovec}\footnote{Mathematical and Computer Sciences and Engineering Division,
		King Abdullah University of Science and Technology,
		Thuwal 23955-6900, Kingdom of Saudi Arabia; 
		{\it jan.haskovec@kaust.edu.sa}}
	}
\vskip 10mm

\noindent{\bf Abstract.}
We rigorously derive the dense graph limit of a discrete model describing the formation of biological transportation networks.
The discrete model, defined on undirected graphs with pressure-driven flows,
incorporates a convex energy functional combining pumping and metabolic costs.
It is constrained by a Kirchhoff law reflecting the local mass conservation.
We first rescale and reformulate the discrete energy functional as an integral `semi-discrete' functional,
where the Kirchhoff law transforms into a nonlocal elliptic integral equation.
Assuming that the sequence of graphs is uniformly connected and that the limiting graphon
is $0-1$ valued, we prove two results: (1) rigorous $\Gamma$-convergence of the sequence
of the semi-discrete functionals to a continuum limit as the number of graph nodes and edges tends to infinity;
 (2) convergence of global minimizers of the discrete functionals
 to a global minimizer of the limiting continuum functional.
Our results provide a rigorous mathematical foundation for the continuum description
of biological transport structures emerging from discrete networks.
\vskip 7mm


%

\vskip 2mm


\vspace{7mm}



\vspace{4mm}
\section{Introduction}\label{sec:Intro}
In this paper we derive the rigorous dense graph limit of the
discrete network formation model introduced and studied in \cite{hu2013adaptation, hu2019optimization}.
The model is posed on a given connected graph $\GraphG = (\Vset,\Eset)$, consisting of a set of vertices (nodes) $\Vset$
and a set of unoriented edges (vessels) $\Eset \subset \Vset\times\Vset$.
Any two vertices $i$, $j\in \Vset$ are linked by at most one edge $(i,j)\in \Eset$, and we also assume that there are no loops, i.e., $(i,i) \notin \Eset$ for all $i \in \Vset$.
Since the graph is undirected we refer by $(i,j)$ and $(j,i)$ to the same edge.
The lengths $\L_{ij} = \L_{ji}>0$ of the edges $(i,j)\in \Eset$
are given and fixed.
Moreover, we denote by $W_{ij}= W_{ji}$ its adjacency matrix, i.e.,
$W_{ij} = 1$ if $(i,j)\in\Eset$ and $W_{ij} = 0$ otherwise.
\par
For each node $i\in \Vset$ we prescribe the strength of source/sink $S_i\in\R$ and we adopt the convention
that $S_i>0$ denotes sources, while $S_i<0$ sinks. We also allow for $S_i=0$, i.e., no external
in- or outgoing flux in this node. We impose the global mass conservation
\(   \label{assS}
   \sum_{i\in \Vset} S_i = 0.
\)
\par
For each vessel $(i,j) \in \Eset$, we denote by $C_{ij}$ its conductivity and by $Q_{ij}$ the flow that runs through it. Since $(i,j)$ and $(j,i)$ refer to the same edge, we have the symmetry $C_{ij} = C_{ji}$ for all $i,j\in \Vset$. Also, $C_{ij} \geq 0$ for all $(i,j) \in \Eset$.
On the other hand, the flux is oriented and antisymmetric, $Q_{ij} = - Q_{ji}$,
and we adopt the convention that $Q_{ij}>0$ means
net flow from node $j\in \Vset$ to node $i\in \Vset$.
\par
We assume low Reynolds number of the flow through the network, so that
the flow rate through a vessel $(i,j)\in \Eset$ is proportional to its conductivity and the pressure drop between its two ends, i.e.,
\(   \label{Q_i}
   Q_{ij} = C_{ij} \frac{P_j - P_i}{\L_{ij}}.
\)
Local conservation of mass is expressed in terms of the Kirchhoff law,
\(  \label{eq:K}
   - \sum_{j\in \mathcal\Vset} W_{ij} C_{ij} \frac{P_j - P_i}{\L_{ij}} = S_i \qquad\mbox{for all } i\in \Vset.
\)
Note that for any given symmetric matrix of nonnegative conductivities $C := (C_{ij})_{(i,j)\in \Vset\times\Vset}$
equation \eqref{eq:K} represents a linear system which is to be solved for the vector of pressures
$P:=(P_i)_{i\in \Vset}$. 
Clearly, a necessary condition for the solvability of \eqref{eq:K} is that the
vector of sources/sinks $S=(S_i)_{i\in\Vset}$ verifies the global mass conservation \eqref{assS}.
Then, the system has a solution for any such $S$ 
if, and only if, the graph with edge weights given by $C$ is connected \cite{Gross-Yellen, HV24},
where only edges with positive conductivities are taken into account (i.e., edges with zero conductivity
are discarded). For a connected graph, the solution is unique up to an additive constant.

The energy cost functional proposed in \cite{hu2013adaptation, hu2019optimization} consists of a pumping power term and a metabolic cost term.
The power (kinetic energy) needed to pump material through an edge $(i,j)\in\Eset$ reads, according to Joule's law and \eqref{Q_i},
\[
    (P_j-P_i) Q_{ij}=\frac{Q_{ij}^2}{C_{ij}}\L_{ij}. 
\]
The metabolic cost of maintaining the edge $(i,j)\in\Eset$ is assumed to be proportional to its length $\L_{ij}$
and a power of its conductivity $C_{ij}^{\gamma}$, with the \emph{metabolic exponent} $\gamma>0$.
In this paper we shall focus on the case $\gamma > 1$.
Consequently, the discrete energy cost functional is given by
\( \label{en:E}
   E[C] := \frac12 \sum_{i\in\Vset} \sum_{j\in\Vset}  \left( C_{ij} \frac{(P_j-P_i)^2}{\L_{ij}^2} + \frac{\nu}{\gamma} C_{ij}^{\gamma} \right)W_{ij} \L_{ij}
\)
where the pressures $P_{ij} = P_{ij}[C]$ are determined as a solution to Kirchhoff's law \eqref{eq:K}, and $\nu>0$ is the so-called metabolic coefficient. If Kirchhoff's law does not admit a solution with the given matrix $C = (C_{ij})_{(i,j)\in \Vset\times\Vset}$,
we set $E[C] := \infty$.

In this paper we derive the rigorous graphon limit of the model \eqref{eq:K}-\eqref{en:E} as the
number of edges $|\Eset|$ and vertices $|\Vset|$ tend to infinity.
The concept of graphons was introduced in \cite{Borgs:2006, Lovasz:2006} to characterize the limiting objects for sequences of dense, finite graphs with respect to the cut metric;
see also the surveys \cite{Lovasz:2012, Glassrock:2015} or the recent paper \cite{ADW:arxiv:2023} for the concept
of probability graphons as limits of large dense weighted graphs. \par
A graphon is a symmetric measurable function $w: [0,1]^2 \to [0,1]$.
In the context of our model \eqref{eq:K}-\eqref{en:E}
we shall interpret $w=w(x,y)$ as the limit of the sequence of graphs
$\GraphG^N = (\Vset^N,\Eset^N)$ as $N\to\infty$, where $|\Vset^N| = N$.
The goal of this paper is then to establish the rigorous limit passage of
the energy functional $E=E[C]$, given by \eqref{eq:K}-\eqref{en:E}, as $N\to\infty$.
In the first step we shall rescale and reformulate the energy functional
in terms of an integral functional acting on the space of piecewise constant
functions on $[0,1]^2$. We then establish its limit in the sense of $\Gamma$-convergence \cite{Maso93}
to an integral functional, where the Kirchhoff law \eqref{eq:K} becomes
an elliptic integral equation in the limit. The `$\limsup$-condition' of $\Gamma$-convergence
follows from strong continuity of the limiting functional with respect to the appropriate topology.
The `$\liminf$-condition' is then a consequence of its weak lower semicontinuity,
which follows from convexity (for values of the metabolic exponent $\gamma\geq 1$). Finally, using the $\Gamma$-convergence,
we shall conclude that global minimizers of the discrete energy functional \eqref{eq:K}-\eqref{en:E}
converge, in appropriate sense, to global minimizers of the limiting (continuum) functional.

To our best knowledge, only three results exist in the literature where continuum
limits of the model \eqref{eq:K}-\eqref{en:E} have been derived.
In \cite{HKM:CMS:2019} a continuum limit is derived formally
for the special case when the graph represents a two-dimensional rectangular grid.
This leads to an integral energy functional on the set of diagonal matrices
(permeability tensors), constrained by a Poisson equation.
In \cite{HKM:CommPDE:2019} this procedure is justified rigorously
in terms of the $\Gamma$-limit of the sequence of properly
rescaled and reformulated energy functionals.
Finally, in \cite{HMP:CMS:2022}, a formal
limit has been derived in the more general setting when
the discrete graphs are triangulations of a bounded two-dimensional domain.
This leads to a similar integral energy functional as in the previous case,
however, defined on the set of symmetric positive semidefinite tensors.
The functional is again constrained by a Poisson equation.

In this paper we derive the graphon limit of \eqref{eq:K}-\eqref{en:E},
where the limiting integral energy functional is defined on the set of nonnegative scalar functions.
These scalar functions represent network conductivities between nodes of the limiting graphon
and are defined on the square $[0,1]^2$, where the interval $[0,1]$
represents node labels. The energy functional is constrained by an integral
Poisson-type equation for the pressure. The gradient of the pressure, which appears
in the ``classical" Poisson equation, is replaced by finite pressure difference between
the corresponding nodes. The scalar conductivity, multiplied by the graphon function,
acts as the integral kernel. Let us stress that our derivation of the limiting functional
as a $\Gamma$-limit of \eqref{eq:K}-\eqref{en:E} is rigorous,
i.e., not only formal, as in \cite{HKM:CMS:2019} and \cite{HMP:CMS:2022}.

This paper is organized as follows.
In Section \ref{sec:notation} we introduce the relevant function spaces
and projection operators and give an overview of the notation that shall be
used throughout the paper.
In Section \ref{sec:ass} we formulate the assumptions on the graph sequence $\GraphG^N$ and its graphon limit.
In Section \ref{sec:results} we reformulate
the discrete functional \eqref{en:E} as an integral one
and state our main results.
In Section \ref{sec:Poisson} we show
that the Kirchhoff law \eqref{eq:K} is equivalent to an integral equation,
which can be seen as a nonlocal version of the Poisson equation.
We establish its well posedness with bounded and square integrable permeability kernels.
In Section \ref{sec:ftionals} we study the integral energy functional
and prove its continuity and convexity.
In Section \ref{sec:thm:1} we provide the proof of our first main result
regarding the $\Gamma$-convergence
of the sequence of reformulated functionals.
Finally, in Section \ref{sec:thm:2} we prove our second main result
about the convergence of the global minimizers.

\vspace{4mm}
\section{Notation and preliminaries}\label{sec:notation}
Throughout the paper we shall denote the real interval $\I := [0,1]$,
and for $N\in\N$ its division into subintervals of length $1/N$,
\[
   \I_i^N := \left[\frac{i-1}{N}, \frac{i}{N} \right]\qquad \mbox{for } i\in [N]
\]
where here and in the sequel we use the notation $[N] := \{1,2,\dots,N\}$.
Moreover, for $N\in\N$ we denote
\[
   \R_0^N := \left\{ v\in \R^N; \, \sum_{i=1}^N v_i = 0 \right\}.
\]
We shall work with square integrable functions on $\I$ with zero mean,
\[
   \Ltwozero(\I) := \left\{ u\in L^2(\I); \, \intI u(x) \d x = 0 \right\},
\]
and with functions in $\Ltwozero(\I)$ that are piecewise constant on each $\I_i^N$,
\[
   \LtwozeroN(\I) := \left\{ u\in L^2_0(\I); \, u \mbox{ is constant on } \I^N_i \mbox{ for all } i\in[N] \right\}.
\]
Moreover, for $q\in [1,\infty]$ and $r > 0$ we introduce the convex set
\(   \label{def:Lqr}
   L^q_r(\II) := \Bigl\{ b\in L^q(\II); \, b(x,y) = b(y,x) \mbox{ and } b(x,y)\geq r \mbox{ for a.e. } x,y\in\I \Bigr\}.
\)

Given a graph $\GraphG^N = (\Vset^N,\Eset^N)$, we
denote by $(W^N_{ij})_{i,j=1}^N$ its adjacency matrix, i.e.,
$W^N_{ij} = 1$ if $(i,j)\in\Eset^N$ and $W^N_{ij} = 0$ otherwise.
For a fixed constant $r>0$ we then define the sequence of sets
\begin{align} \label{def:BNr}
\BB^N_r := \left\{ B \in \R^{N \times N}_{\mathrm{sym}} \,\middle|\,
\begin{array}{l}
B_{ij} \geq r \quad \text{for all } (i,j) \in [N]^2 \text{ such that } W^N_{ij} = 1, \\
B_{ij} = 0 \quad \text{for all } (i,j) \in [N]^2 \text{ such that } W^N_{ij} = 0
\end{array}
\right\}.
\end{align}
Obviously, the sets $\BB^N_r$ are convex for all $N\in\N$.

For any function $p$ defined on $\I$ we shall use the notation
\(   \label{def:notation:dd}
    \dd p(x,y) := p(y)-p(x) \qquad\mbox{for } (x,y)\in \II.
\)

For $N\in\N$ we introduce the operators $\QQ^N: \R^{N} \to L^\infty(\I)$
mapping vectors $v = (v_1,\ldots,v_N)\in\R^N$ onto piecewise constant functions on $\I$,
\[     
   \QQ^N[v](x) := \sum_{i=1}^N v_i \, \chi^N_i(x)\qquad \mbox{for } x\in \I,
\]
where $ \chi^N_i = \chi_{\I^N_i}$ is the characteristic function of the interval $\I^N_i$.
By a slight abuse of notation, we also use the symbol $\QQ^N$ for the operator
$\QQ^N: \R^{N\times N} \to L^\infty(\II)$, mapping matrices $B\in\R^{N\times N}$
onto piecewise constant functions on $\II$,
\[   
   \QQ^N[B](x,y) := \sum_{i=1}^N\sum_{j=1}^N B_{ij} \, \chi^N_i(x)\chi^N_j(y) \qquad \mbox{for } (x,y)\in \I\times\I.
\]
For each adjacency matrix $W^N \in \R^{N\times N}$
we introduce the so-called pixel picture,
which is the bounded piecewise constant function
$w^N := \QQ^N[W^N]$.
Moreover, throughout the paper we use the notation
\[
   \ell^N := \QQ^N[\L^N].
\]

We also define the operator $\ZZ^N : L^1(\I) \to L^\infty(\I)$,
being a projection on the space of piecewise constant functions
on the intervals $\I_i^N$, i.e., for any function $u\in L^1(\I)$,
\[   
     \ZZ^N[u](x) := N \int_{\I_i^N} u(s) \d s \qquad\mbox{for } x \in \I_i^N.
\]
We note that $\ZZ^N$ is an orthogonal projection on $L^2(\I)$
and as such is self adjoint, i.e., for any $u$, $v\in L^2(\I)$ we have
\(   \label{eq:ZZ:ortho}
   \intI \ZZ^N[u] v \d x = \intI u \ZZ^N[v] \d x.
\)
Again, by an abuse of notation, we use the same symbol for the operator
$\ZZ^N : L^1(\II) \to L^\infty(\II)$,
mapping onto the space of piecewise constant functions
on the patches $\I_i^N \times \I_j^N$, i.e., for any function $u\in L^1(\II)$,
\[   
     \ZZ^N[u](x,y) := N^2 \int_{\I_i^N} \int_{\I_j^N} u(s,t) \d s\d t \qquad\mbox{for } (x,y) \in \I_i^N\times\I_j^N.
\]
The self-adjoint property \eqref{eq:ZZ:ortho} also holds in this case with an obvious modification.

\vspace{4mm}
\section{Main assumptions}\label{sec:ass}
We now give the overview of assumptions
that we impose on the sequence of graphs
$\GraphG^N = (\Vset^N,\Eset^N)$ and the source/sink terms $S^N\in \R^N_0$.
Let us recall that we denote $W^N\in \R^{N\times N}$ the adjacency matrix of $\GraphG^N$
and $w^N := \QQ^N[W^N]$ its pixel picture.

We start by imposing the assumption that there exists a function $\sigma\in \Ltwozero(\I)$ such that
the source/sink terms are given by 
\( \label{ass:S}
   S_i^N := \int_{\I_i^N} \sigma(x) \d x \qquad\mbox{for } i\in[N].   
\)

Next, we introduce two assumptions on the adjacency matrices and their pixel pictures. 
Namely, we assume that 
    there exists $\lambda>0$ such that for all $N\in\N$,
    \(  \label{ass:A1}
        \sum_{i=1}^N \sum_{j=1}^N  (z_i - z_j)^2 W^N_{ij} \geq \lambda N \sum_{i=1}^N z_i^2  \qquad\mbox{for all } z\in\R_0^N.
    \)
Moreover, we assume that
    there exists $w\in L^\infty(\II)$ such that
    \(  \label{ass:A2}
       w^N \to w \mbox{ in the norm topology of } L^1(\II) \mbox{ as } N\to\infty.
    \)
    
Assumption \eqref{ass:A1} is equivalent to the lower bound $\mathfrak{f}(\GraphG^N) \geq \lambda N$
on the Fiedler number $\mathfrak{f}(\GraphG^N)$ of $\GraphG^N$.
The Fiedler number, also called algebraic connectivity,
is the second smallest eigenvalue of the matrix Laplacian of $W^N$; see, e.g., \cite{Fiedler, Mieghem}.
Therefore, we can say that \eqref{ass:A1} enforces uniform connectivity of the graphs $\GraphG^N$.
A particular example follows from \cite[Proposition 3.8]{Fiedler}, which states that
\[  
   \mathfrak{f}(\GraphG^N) \geq 2 \delta(\GraphG^N) - N + 2,
\]
where $\delta(\GraphG^N)$ is the minimal vertex degree of $\GraphG^N$.
Obviously, \eqref{ass:A1} is verified if the minimal vertex degree
of all the graphs $\GraphG^N$, $N\in\N$, is at least $\frac{(1+\lambda)N}{2}$.
Then we have
\[
   \mathfrak{f}(\GraphG^N) \geq 2 \frac{(1+\lambda)N}{2}  - N + 2 = \lambda N + 2.
\]

Assumption \eqref{ass:A2} restricts the validity of our results
to $0-1$ valued graphons, i.e., $\mbox{range}(w) \subset \{0, 1\}$.
In this case one can show \cite[Proposition 8.24]{Lovasz:2012} that the usual convergence of $w^N$ to $w$
in the cut norm indeed implies convergence in the $L^1$ norm topology.
Of course, due to the uniform boundedness of $w^N$ in $L^\infty(\II)$,
the convergence \eqref{ass:A2} holds also in the norm topology of $L^q(\II)$ for any $q<\infty$.
Let us explicitly note that \eqref{ass:A2} excludes important classes of graphs.
For instance, for the generic random Erd\"{o}s-R\'{e}nyi graph model \cite{Erdos-Renyi},
where pair of nodes are connected with probability $p\in (0,1)$,
we have $w\equiv p$ on $\II$. Obviously, it is then not possible that $w^N \to w$ strongly in $L^1$,
since the functions $w^N$ are piecewise constant taking only values $0$ and $1$.
We shall generalize our theory to include random graph models in a future work.

To be able to pass to the continuum limit in \eqref{eq:K}-\eqref{en:E},
we also need to assume that the edge lengths $\L^N =(\L^N_{ij})_{i,j=1}^N$ converge
in an appropriate sense.
In particular, for $\ell^N:=\QQ^N[\L^N]$ we impose the assumption that
\(  \label{ass:L1}
   \Norm{\ell^N}_{L^\infty(\II)} = \sup_{i,j\in[N]} \L^N_{ij} \leq 1,
\)
and, moreover, there exists $\ell\in L^\infty(\II)$ such that
\(  \label{ass:L2}
   \ell^N \to \ell \mbox{ in the norm topology of } L^1(\II) \mbox{ as } N\to\infty.
\)
Assumption \eqref{ass:L1} is a natural assumption meaning that the sequence of
graphs $\GraphG^N$ is located within a compact set in the ``physical space".
The particular bound $\L_{ij} \leq 1$ can be achieved through rescaling and is therefore
assumed without loss of generality.
A trivial way to ensure validity of assumption \eqref{ass:L2} is to
fix $\ell\in L^\infty(\II)$ and set $\ell^N:=\ZZ^N[\ell]$.
The lengths $\L^N_{ij}$ are then the corresponding values of $\ell^N$,
i.e., such that $\ell^N:=\QQ^N[\L^N]$ holds.
Again, due to the uniform boundedness of $\ell^N$ in $L^\infty(\II)$,
the limit in \eqref{ass:L2} holds also with respect to the norm topology of $L^q(\II)$ for any $q<\infty$.

Assumptions \eqref{ass:A1}--\eqref{ass:L2} are sufficient for proving our first main result,
which establishes the $\Gamma$-convergence of the sequence of (properly rescaled and reformulated)
energy functionals. For our second main result, which states that the sequence of global minimizers
of the discrete functionals \eqref{eq:K}-\eqref{en:E} converges to the global minimizer of the limiting functional,
we need to adopt one more assumption. Namely, we need to assume that $\gamma>2$ and that
the sequence of reciprocal values $\left(\ell^N\right)^{-1}$
is uniformly bounded in $L^{\frac{2(\gamma+1)}{\gamma-2}}(\II)$,
\begin{equation}\label{ass:L3star}
   \sup_{N\in\N} \Norm{\left(\ell^N\right)^{-1}}_{L^{\frac{2(\gamma+1)}{\gamma-2}}(\II)} < \infty.
\end{equation}
Admittedly, this is a rather strong assumption and constitutes a key limitation of our approach.
To gain an intuition about the satisfiability of \eqref{ass:L3star},
let us consider the case when the vertices in $\Vset^N$
are generated as random points in a convex compact set $K\subset \R^d$.
I.e., we identify $\Vset^N$ with a cloud of $N\in\N$ points drawn independently from the
uniform distribution on $K$.
The edge lengths $\L_{ij}$ are then simply
the Euclidean distances between the points (vertices) $i$ and $j$.
The law of large numbers~\cite{Billingsley:95} gives then
\begin{equation}\label{eq:LLN}
       \lim_{N\to\infty} \frac{1}{N^2} \sum_{i=1}^N \sum_{j\neq i}^N \L_{ij}^{-q}
       = \int_0^{\infty} s^{-q} f(s) \d s,
\end{equation}
where $f=f(s)$ is the probability density for the Euclidean distances of the point cloud.
For the case when $K$ is a hypercube in $\R^d$,
the probability density $f=f(s)$ is known explicitly; see, e.g., \cite{Kendall:1960, LS:Stats:2020} and references therein.
For $d=2$ and $K = [0,1/\sqrt{2}]^2$, we have $\L_{ij}\leq 1$ and
\[
   f(s) = 8s \left(s^2-2\sqrt{2}s+\frac{\pi}{2}\right) \qquad\mbox{for } s\in\left[0,\frac{1}{\sqrt{2}}\right],
\]
complemented by a more complicated formula for values of $s\in \left(\frac{1}{\sqrt{2}}, 1\right]$.
However, $f(s)$ is decreasing for $s>1/\sqrt{2}$ (see~\cite{LS:Stats:2020}) so that
we have the bound $f(s) < f\left(1/\sqrt{2}\right) = 2\sqrt{2}(\pi-3)$ for $s>1/\sqrt{2}$.
Of course, $f(s)=0$ for all $s>1$.
Then, in \eqref{eq:LLN} we have
\[
    \int_0^\infty s^{-q} f(s) \d s 
     \leq  8 \int_0^{1/\sqrt{2}} s^{-q+1} \left(s^2-4\sqrt{2}s+\frac{\pi}{2}\right) \d s
     + 2\sqrt{2}(\pi-3) \int_{1/\sqrt{2}} ^1 s^{-q} \d s.
\]
Obviously, the right-hand side is finite if and only if the first integral is finite, which holds if and only if $q<2$.
Unfortunately, this makes assumption \eqref{ass:L3star} unsatisfiable in dimension $d=2$,
since
\[
   \frac{2(\gamma+1)}{\gamma-2} > 2
   \qquad\mbox{for all }\gamma>2.
\]
In general, the probability density for the $d$-dimensional cube
behaves like $s^{d-1}$ as $s\to 0+$; see~\cite{LS:Stats:2020}.
Consequently, \eqref{eq:LLN} is finite if and only if $q<d$.
Therefore, assumption \eqref{ass:L3star} requires
\[
   \frac{2(\gamma+1)}{\gamma-2} < d,
\]
which is equivalent to
\[
   \gamma > \frac{2d+2}{d-2} = 2 + \frac{6}{d-2}.
\]
In particular, for $d=3$ we need $\gamma>8$,
while for $d=4$ we need $\gamma>5$, and for $d=5$ we need $\gamma>4$.


\vspace{4mm}
\section{Reformulation of the model and main results}\label{sec:results}
For the purpose of deriving the graphon limit
of the model \eqref{eq:K}-\eqref{en:E}, we first introduce a regularization
by adding a fixed constant $r>0$ to all conductivities $C_{ij}$.
This is necessary to ensure uniform ellipticity of the Kirchhoff law \eqref{eq:K}
and was also introduced in the previous works
\cite{HKM:CMS:2019, HKM:CommPDE:2019} and \cite{HMP:CMS:2022}
where continuum limits were derived.
It is then convenient to introduce the new set of variables $B_{ij} := \frac{C_{ij} + r}{\L_{ij}}$.
Moreover, the energy functional $E^N$ has to be scaled
by $1/N^2$, which can be interpreted as the reciprocal of the
square of the ``mean" edge length.
This is due to the fact that we are embedding the graph $\GraphG^N$,
which is inherently a one-dimensional structure,
into the two-dimensional interval $\II$.
The same scaling was adopted in \cite{HKM:CMS:2019, HKM:CommPDE:2019,HMP:CMS:2022} in order to obtain a meaningful continuum limit.
We therefore introduce the rescaled energy functional
\(  \label{eq:F}
   F^N[B] := \frac1{2N^2} \sum_{i=1}^N \sum_{j=1}^N \left( B_{ij} (P_j-P_i)^2 + \frac{\nu}{\gamma} B_{ij}^{\gamma}  \left( \L_{ij}^N \right)^{\gamma+1} \right) W^N_{ij},
\)
coupled to the rescaled Kirchhoff law
\(   \label{eq:KB}
   - \frac{1}{N^2} \sum_{j=1}^N W^N_{ij} B_{ij} (P_j - P_i) = S_i^N \quad\mbox{for all } i\in [N].
\)
Clearly, we have $F^N[B] = E^N[C+r]/N^2$, with $(C+r)_{ij}:=C_{ij}+r$.
We fix $r>0$ and, recalling the assumption $\L_{ij}\leq 1$, pose \eqref{eq:F}--\eqref{eq:KB} on the convex set $\BB^N_r$ defined in \eqref{def:BNr}.
We shall also use the ``weak formulation" of \eqref{eq:KB},
which is obtained by its multiplication by $\Phi_i$ for a ``test vector" $\Phi\in\R_0^N$ and summation over $i\in[N]$,
\(   \label{eq:KBw}
   \frac1{2N^2} \sum_{i=1}^N \sum_{j=1}^N W_{ij}^N B_{ij} (P_j - P_i) (\Phi_j - \Phi_i) = \sum_{i=1}^N S_i^N \Phi_i.
\)

In the continuum limit we shall derive the integral functional
\(  \label{eq:Fcont}
   \Fcont[\b] =  \frac12 \iintII  b(x,y) \left(p(x)-p(y)\right)^2 + \frac{\nu}{\gamma} b(x,y)^\gamma \ell(x,y)^{\gamma+1}  \d w(x,y),
\)
where we use the concise notation $\d w(x,y)$ for $w(x,y)\d x\d y$ and $w\in L^\infty(\II)$ is given in \eqref{ass:A2}.
The functional $\Fcont$ is defined on the set $L^\infty_r(\II)$ introduced in \eqref{def:Lqr}.
The scalar pressure $p=p[\b]\in \Ltwozero(\I)$ is obtained as a solution of the nonlocal diffusion (Poisson-type) equation
\[  
   \intI \b(x,y) (p(x)-p(y))  w(x,y)\d y = \sigma(x),
\]
where the datum $\sigma\in \Ltwozero(\I)$ represents the intensity of sources and sinks;
cf. \cite{Medvedev2014, Medvedev:SIAM, EPSS:ARMA:2021}.
We shall work with its symmetrized weak formulation, which is obtained
by a multiplication by a test function $\varphi\in \Ltwozero(\I)$
and integration over $\I$,
\( \label{eq:Poissonbw}
   \frac12 \iintII \b(x,y) (\dd p(x,y)) (\dd \varphi(x,y)) \d w(x,y) = \int_\I \sigma(x) \varphi(x) \d x \quad \forall \varphi\in \Ltwozero(\I),
\)
where we used the notation \eqref{def:notation:dd}.
Let us note that the functional \eqref{eq:Fcont}--\eqref{eq:Poissonbw} is invariant
with respect to any measure preserving change of variable on the interval $\I$.
Indeed, if $\psi:\I\to\I$ is any such measure preserving bijection,
then the value of the energy remains the same for the transformed functions
$\tilde b(x,y) := b(\psi(x),\psi(y))$, $\tilde\ell(x,y) := \ell(\psi(x),\psi(y))$, $\tilde w(x,y):=w(\psi(x),\psi(y))$
and $\tilde p(x) := p(\psi(x))$, $\tilde\sigma(x):=\sigma(\psi(x))$.
This invariance is important since $w$ and $\tilde w$ represent the same graphon \cite{Borgs:2006, Lovasz:2012}.
An analogous invariance obviously also holds for the discrete model \eqref{eq:F}--\eqref{eq:KB} with
respect to permutations of the set of indices $[N]$.

To establish a connection between the discrete functional \eqref{eq:F}--\eqref{eq:KB} and
the continuum limit \eqref{eq:Fcont}--\eqref{eq:Poissonbw},
we introduce for $N\in\N$ the intermediary ``semi-discrete" energy functional
\(  \label{Fn:N}
   \Fcont^N[b] :=  \frac12 \intI\intI b(x,y) \left(\dd p^N(x,y)\right)^2 
         + \frac{\nu}{\gamma} b(x,y)^\gamma \ell^N(x,y)^{\gamma+1} \d w^N(x,y),
\)
with $\ell^N=\QQ^N[\L^N]$, $w^N=\QQ^N[W^N]$,
and $p^N\in\LtwozeroN(\I)$ is the unique piecewise constant solution of the approximate Poisson equation
\begin{equation}\label{eq:Poissonsc}
\frac12 \iintII \b(x,y) (\dd p^N(x,y)) (\dd \varphi^N(x,y)) \d w^N(x,y) = \int_\I \sigma(x) \varphi^N(x) \d x \quad \forall \varphi^N\in \LtwozeroN(\I).
\end{equation}
\par
We are now ready to formulate our main results.

\begin{theorem}\label{thm:1}
Fix $\gamma>1$, $\nu > 0$, $\lambda>0$, and $r>0$, and let the assumptions \eqref{ass:S}--\eqref{ass:L2} hold.
For any $N\in\N$ and $B\in \BB^N_r$ we have
\(   \label{eq:connection}
   \Fcont^N[\QQ^N[B]] = F^N[B].
\)
Moreover, the sequence of functionals $(\Fcont^N)_{N\in\N}$, given by \eqref{Fn:N}--\eqref{eq:Poissonsc},
$\Gamma$-converges in $L_r^\omega(\II)$, $\omega=\max\{2,\gamma\}$,
to the functional $\Fcont$ given by \eqref{eq:Fcont}--\eqref{eq:Poissonbw},
in the following sense:
\begin{enumerate}[label=(\roman*)]
\item
For any sequence $(\b^N)_{N \in \N}\subset L_r^\omega(\II)$ converging weakly 
in $L^\omega(\II)$ to $\b\in L_r^\omega(\II)$ we have
\(   \label{eq:gamma:liminf}
   \Fcont[\b] \leq \liminf_{N\to\infty}  \Fcont^N[\b^N].
\)
\item
For any $\b\in L_r^\omega(\II)$ the sequence $\b^N := \ZZ^N[\b]$ 
converges to $\b$ in the norm topology of $L^\omega(\II)$,
and
\(   \label{eq:gamma:limsup}
   \Fcont[\b] = \lim_{N\to\infty}  \Fcont^N[\b^N].
\)
\end{enumerate}
\end{theorem}

The first claim of Theorem \ref{thm:1}, formula \eqref{eq:connection},
establishes a connection between the discrete energy functional $F^N$
and its ``semi-discrete" counterpart $\Fcont^N$. The second claim then
characterizes the convergence of $\Fcont^N$ as $N\to\infty$.

Our second main result establishes the convergence of the sequence
of global minimizers of the discrete functional $F^N$
to the global minimizer of the continuum energy functional $\Fcont$.

\begin{theorem}\label{thm:2}
Fix $\gamma>2$, $\nu > 0$, $\lambda>0$, and $r>0$, and let the assumptions \eqref{ass:S}--\eqref{ass:L3star} hold.
Let $B^N \in \BB^N_r$, $N\in\N$, be a sequence of global minimizers of the discrete energy functionals $F^N$,
given by \eqref{eq:F}--\eqref{eq:KB}, on the sets $\BB^N_r$.
Then $\b^N = \QQ^N[B^N]$ converges weakly in $L^2(\II)$
to $\b\in L^2(\II)$, which is the global minimizer of the continuum energy functional
$\Fcont$, given by \eqref{eq:Fcont}--\eqref{eq:Poissonbw}, on $L_r^\gamma(\II)$.
\end{theorem}

We note that unique global minimizers of the discrete energy functionals $F^N$
do exist for each $N\in\N$ due to their continuity, strict convexity (with $\gamma>1$) and coercivity
on the convex sets $\BB^N_r$; see \cite[Appendix]{HV24} for details.
Moreover, the global minimizers of the continuum energy functional
$\Fcont$ can be found using the corresponding gradient flow.
Indeed, analogously to \cite[Lemma 3.1]{HKM:CMS:2019},
one can formally derive the following explicit formula for the $L^2$-gradient flow
of \eqref{eq:Fcont},  constrained by \eqref{eq:Poissonbw},
\[   
  \partial_t b(x,y) = (p(x)-p(y))^2 - \nu b(x,y)^{\gamma-1} \ell(x,y)^{\gamma+1}.
\]

\vspace{4mm}
\section{The Kirchhoff law and the Poisson equation}\label{sec:Poisson}

We start by formulating the following simple Poincar\'{e}-type equality for functions with vanishing mean over $\I$.

\begin{lemma} \label{lem:Poincare}
For any $p\in \Ltwozero(\I)$ we have
\(  \label{eq:Poincare}
   \sqrt{2} \Norm{p}_{L^2(\I)} = \Norm{\dd p}_{L^2(\II)}.
\)
\end{lemma}

\begin{proof}
With $\intI p(y) \d y = 0$ and $\mbox{meas}(\I)=1$, we have
\[
   \Norm{\dd p}_{L^2(\II)}^2 &=& \iintII (p(x) - p(y))^2 \d x \d y \\
     &=&  \iintII p(x)^2 - 2p(x)p(y) + p(y)^2 \d x\d y = 2 \Norm{p}^2_{L^2(\I)}.
\]
\end{proof}

Next we derive a simple estimate on the solutions of the Kirchhoff law \eqref{eq:KB}.

\begin{lemma}\label{lem:Kirchhoff}
Let $\lambda>0$, $r>0$ and $S^N\in\R^N_0$ be given by \eqref{ass:S} with some $\sigma\in \Ltwozero(\I)$.
Then for any conductivity matrix $B^N \in \BB^N_r$ there exists
a unique solution $P\in \R^N_0$ of the Kirchhoff law \eqref{eq:KB}.
Moreover, we have the estimate
\(  \label{est:Kirchhoff}
   \sum_{i=1}^N \sum_{j=1}^N (P_i - P_j)^2 \leq \frac{8N^2}{(r\lambda)^2} \intI \sigma(x)^2 \d x.
\)
\end{lemma}

\begin{proof}
The unique solvability of \eqref{eq:KB} in $\R^N_0$
is a direct consequence of the Lax-Milgram theorem,
where the coercivity follows from assumption \eqref{ass:A1}
and the fact that $B_{ij} \geq r$ for all $i,j\in[N]$.

We now use the solution $P\in \R^N_0$ as the test vector $\Phi$ in the symmetrized formulation \eqref{eq:KBw},
which gives
\[
   \frac1{2N^2} \sum_{i=1}^N \sum_{j=1}^N W_{ij}^N B^N_{ij} (P_j - P_i)^2 = \sum_{i=1}^N S^N_i P_i.
\]
Using again \eqref{ass:A1}, the left-hand side is estimated from below by
\[
   \frac1{2N^2} \sum_{i=1}^N \sum_{j=1}^N W_{ij}^N B^N_{ij} (P_j - P_i)^2 \geq \frac{r \lambda}{2N} \sum_{i=1}^N P_i^2.
\]
The right-hand side is estimated using the Cauchy-Schwarz inequality
\[
   \sum_{i=1}^N S^N_i P_i \leq \left( \sum_{i=1}^N \left(S^N_i\right)^2 \right)^\frac12 \left( \sum_{i=1}^N P_i^2 \right)^\frac12.
\]
Consequently, we have
\[
   \sum_{i=1}^N P_i^2 \leq \left( \frac{2N}{r\lambda}\right)^2 \sum_{i=1}^N \left(S^N_i\right)^2.
\]
Since $P\in \R^N_0$, we have
\[
   \sum_{i=1}^N \sum_{j=1}^N (P_i - P_j)^2 = 2 N \sum_{i=1}^N P_i^2,
\]
and \eqref{ass:S} gives
\[
   \sum_{i=1}^N \left(S^N_i\right)^2 = \sum_{i=1}^N \left( \int_{\I_i^N} \sigma(x) \d x \right)^2 \leq \frac{1}{N} \intI \sigma(x)^2 \d x.
\]
A concatenation of the above estimates directly gives \eqref{est:Kirchhoff}.
\end{proof}

We  have the following simple consequence of assumptions \eqref{ass:A1} and \eqref{ass:A2}.

\begin{lemma}\label{lem:w}
Let assumptions \eqref{ass:A1} and \eqref{ass:A2} be verified.
Then
\(  \label{ass:A3}
    \iintII  \big(\dd z(x,y)\big)^2  \d w(x, y) \geq \lambda \intI z(x)^2 \d x \qquad\mbox{for all } z\in L^2_0(\I).
\)
\end{lemma}

\begin{proof}
Let us fix $z\in L^2_0(\I)$.
For any $N\in\N$, denoting $\ZZ^N[z]_i$ the constant value of the function $\ZZ^N[z]$ on $\I_i^N$,
we readily have
\[
   \iintII  \left(\ZZ^N[z](x) - \ZZ^N[z](y) \right)^2  \d w^N(x, y) = \frac{1}{N^2} \sum_{i=1}^N \sum_{j=1}^N \left( \ZZ^N[z]_i - \ZZ^N[z]_j \right)^2 W_{ij}^N,
\]
and
\[
   \intI \ZZ^N[z](x)^2 \d x = \frac{1}{N} \sum_{i=1}^N \left( \ZZ^N[z]_i \right)^2.
\]
Assumption \eqref{ass:A1} gives then
\[
   \iintII  \left(\ZZ^N[z](x) - \ZZ^N[z](y) \right)^2  \d w^N(x, y) \geq \lambda \intI \ZZ^N[z](x)^2 \d x.
\]
The claim \eqref{ass:A3} is obtained by passing to the limit $N\to \infty$.
For this we use the strong convergence of $\ZZ^N[z]$ to $z$ in $L^2(\I)$,
see, e.g., \cite{Rudin},
assumption \eqref{ass:A2} and the uniform boundeness of $w^N$ in $L^\infty(\II)$.
\end{proof}

Now we establish a connection between solutions of the Kirchhoff law \eqref{eq:KB}
and weak solutions of the integral equations \eqref{eq:Poissonbw} and
\eqref{eq:Poissonsc}. 

\begin{lemma}\label{lem:Poisson}
Let $\sigma\in \Ltwozero(\I)$ and $\b\in L^\infty_r(\II)$ for some $r>0$. Assume \eqref{ass:A1} and \eqref{ass:A2}.
Then there exists a unique $p\in \Ltwozero(\I)$ verifying \eqref{eq:Poissonbw},
and a unique $p^N \in \LtwozeroN(\I)$ verifying \eqref{eq:Poissonsc}.

Moreover, if $\b = \QQ^N[B]$ for some $B\in\BB_r^N$,
then we have $p^N = \ZZ^N[p]$. Denoting
\(  \label{P_i}
   P_i := p^N(x) = \ZZ^N[p](x) \qquad\mbox{for } x\in\I_i^N,\; i\in[N],
\)
then $(P_i)_{i=1}^N \in\R_0^N$ is a solution of Kirchhoff's law \eqref{eq:KB} with conductivities $(B_{i,j})_{i,j=1}^N$
and source/sink terms $S^N\in \R_0^N$ given by \eqref{ass:S}.
\end{lemma}

\begin{proof}
The existence and uniqueness of the solutions of \eqref{eq:Poissonbw}
and \eqref{eq:Poissonsc} follows directly from the Lax-Milgram theorem.
The coercivity of the respective bilinear forms on $L^2_0(\I)$ follows from \eqref{ass:A3} together with the lower bound $b\geq r>0$ almost everywhere on $\II$ and the Poincar\'{e}-type equality \eqref{eq:Poincare}.

In case of the piecewise constant $\b = \QQ^N[B]$, we have $\b(x,y) \equiv B_{ij}$ for all $(x,y)\in\I_i^N\times\I_j^N$.
Let $p\in \Ltwozero(\I)$ be the unique solution of \eqref{eq:Poissonbw}. Taking a piecewise constant test function
$\varphi\in \LtwozeroN(\I)$, we have for the left-hand side of \eqref{eq:Poissonbw},
\begin{align*}
 \frac12 \iintII  \b(x,y) & (\dd p(x,y)) (\dd \varphi(x,y)) \d w^N(x, y) \\ & = 
   \frac12 \sum_{i=1}^N \sum_{j=1}^N (\varphi_j-\varphi_i) W_{ij}^N B_{ij} \int_{\I_i^N} \int_{\I_j^N} \dd p(x,y) \d x\d y,
\end{align*}
where we denoted $\varphi_i$ the constant value of $\varphi$ on $\I_i^N$, and similarly for $\varphi_j$.
Moreover,
\[
   \int_{\I_i^N}\int_{\I_j^N} \dd p(x,y) \d x\d y = \frac{1}{N} \int_{\I_j^N} p(y) \d y - \frac{1}{N} \int_{\I_i^N} p(x) \d x = \frac{1}{N^2} \left( \ZZ^N[p]_j - \ZZ^N[p]_i \right),
\]
where we again denoted $\ZZ^N[p]_i$ the constant value of $\ZZ^N[p]$ on $\I_i^N$.
For the right-hand side of \eqref{eq:Poissonbw} we have
\(  \label{eq:PoissonRHS}
    \intI \sigma(x) \varphi(x) \d x = \sum_{i=1}^N \varphi_i \int_{\I_i^N} \sigma(x) \d x = \sum_{i=1}^N \varphi_i S_i^N.
\)
We conclude
\[
   \frac1{2N^2} \sum_{i=1}^N \sum_{j=1}^N  W_{ij}^N  B_{ij} \left( \ZZ^N[p]_i - \ZZ^N[p]_j \right) (\varphi_i-\varphi_j) = \sum_{i=1}^N \varphi_i S_i^N,
\]
which is the weak formulation \eqref{eq:KBw} with the test vector $\varphi$,
and thus $\left(\ZZ^N[p]_i\right)_{i=1}^N$ is a solution of the Kirchhoff law \eqref{eq:KB}.

Finally, if $p^N \in \LtwozeroN(\I)$ is the unique solution of \eqref{eq:Poissonsc}
and $\varphi\in \LtwozeroN(\I)$ is a piecewise constant test function,
then, trivially,
\begin{align*}
   \frac12 \iintII  \b(x,y) & (\dd p^N(x,y)) (\dd \varphi(x,y)) \d w^N(x, y) \\
   &  =
   \frac1{2N^2} \sum_{i=1}^N \sum_{j=1}^N (\varphi_i-\varphi_j) (p^N_i - p^N_j) W_{ij}^N  B_{ij},
\end{align*}
where we denoted $p^N_i$ the constant value of $p^N$ on $\I_i^N$, and similarly for $p_j^N$.
The right-hand side of \eqref{eq:Poissonsc} is identical to \eqref{eq:PoissonRHS}.
We conclude that $\left(p^N_i\right)_{i=1}^N$ is a solution of the Kirchhoff law \eqref{eq:KB}, too.
Due to the uniqueness of the solutions of the Kirchhoff law (up to an additive constant),
and since $\int_\I p(x) \d x = \int_\I p^N(x) \d x = 0$,
we have $p^N = \ZZ^N[p]$ and \eqref{P_i} follows.
\end{proof}

Finally, we establish solvability of the Poisson equation in its weak formulation
\eqref{eq:Poissonbw} with a square-integrable permeability kernel $\b\in L^2_r(\II)$.

\begin{lemma}\label{lem:Poisson-b}
Let $\sigma \in L^2_0(\I)$ and $\b \in L^2_r(\II)$. Assume \eqref{ass:A1} and \eqref{ass:A2}.
Then there exists a unique $p \in \Ltwozero(\I)$ verifying \eqref{eq:Poissonbw}
for all test functions $\varphi\in L_0^\infty(\I)$.
Moreover, we have
\(  \label{eq:testPoisson}
   \frac12 \iintII \b (\dd p)^2 \d w(x,y) = \intI \sigma p \d x < \infty
\)
and
\(   \label{est:p}
     \Norm{p}_{L^2(\I)}  \leq  \frac{2}{\lambda r} \Norm{\sigma}_{L^2(\I)}.
\)
\end{lemma}

\begin{proof}
For $n\in\N$ we define the cut-off kernel $\b^n(x,y) := \min \{\b(x,y), n \}$, $(x,y)\in\II$.
Then $\b^n\in L^\infty_r(\II)$ for all $n\in\N$, $n>r$, and Lemma \ref{lem:Poisson} provides a unique $p^n \in \Ltwozero(\I)$
which verifies \eqref{eq:Poissonbw} with kernel $b^n$ and test functions $\varphi\in \Ltwozero(\I)$.

Using $p^n$ as a test function in \eqref{eq:Poissonbw}, due to \eqref{ass:A3} and the uniform boundedness of $\b^n$ from below,
we have
\(   \label{auxEst2}
  \lambda r \Norm{p^n}_{L^2(\I)}^2 \leq \iintII \b^n |\dd p^n|^2 \d w(x,y) &=& 2\int_\I \sigma p^n \d x \\
      &\leq& 2 \Norm{\sigma}_{L^2(\I)} \Norm{p^n}_{L^2(\I)},  \nonumber
\)
We thus have the uniform bound 
\(   \label{est:pn}
   \Norm{p^n}_{L^2(\I)}   \leq  \frac{2}{\lambda r} \Norm{\sigma}_{L^2(\I)}
\)
and there exists a subsequence of $p^n$ that converges weakly in $L^2(\I)$ to some $p\in \Ltwozero(\I)$.
By the Lebesgue dominated convergence theorem, $\b^n\to b$ strongly in $L^2(\II)$ as $n\to \infty$.
We can therefore pass to the limit in the weak formulation \eqref{eq:Poissonbw}
to obtain
\(  \label{p-distr2}
   \frac12 \intI \intI  \b(x,y) \dd p(x,y) \dd\varphi(x,y) \d w(x,y) = \intI \sigma(x) \varphi(x) \d x
\)
for all test functions $\varphi\in L_0^\infty(\I)$.
Noting that \eqref{auxEst2} also implies a uniform bound on $\iintII \b^n |\dd p^n|^2 \d w(x,y)$,
we have due to the weak lower semicontinuity of the $L^2$-norm,
\[
   \iintII \b (\dd p)^2 \d w(x,y) = \Norm{ \sqrt{\b w} (\dd p) }_{L^2(\II)}^2
    \leq \liminf_{n\to\infty} \Norm{ \sqrt{\b^n w} (\dd p^n) }_{L^2(\II)}^2
    < \infty,
\]
where for the identification of the limit we used the fact that $\sqrt{\b^n} \to \sqrt{\b}$
in the norm topology of $L^4(\II)$, recalling that $b^n \geq r >0$,
and $\dd p^n \wto \dd p$ weakly in $L^2(\II)$.
Consequently, we can use a sequence of cut-off versions of $p$ as test functions in \eqref{p-distr2} and pass to the limit
to remove the cut-off. This gives
\[
   \frac12 \iintII b (\dd p)^2 \d w(x,y) = \intI \sigma p \d x,
\]
which is \eqref{eq:testPoisson}.
Estimate \eqref{est:p} follows from \eqref{est:pn} and the weak lower semicontinuity of the $L^2$ norm.
\end{proof}

An obvious modification of the above proof yields
the following result for the approximate Poisson equation \eqref{eq:Poissonsc}.

\begin{lemma}\label{lem:Poisson-bh}
Let $\sigma \in L^2_0(\I)$ and $\b \in L^2_r(\II)$. Assume \eqref{ass:A1} and \eqref{ass:A2}.
Then for any $N\in\N$ there exists a unique $p^N \in \LtwozeroN(\I)$ verifying \eqref{eq:Poissonsc}
for all test functions $\varphi^N\in L^\infty_{0,N}(\I)$.
Moreover, we have
\(  \label{eq:testPoisson-h}
   \frac12 \iintII \b (\dd p^N)^2 \d w^N(x,y) = \intI \sigma p^N \d x < \infty,
\)
and
\[   
     \Norm{p^N}_{L^2(\I)}  \leq  \frac{2}{\lambda r} \Norm{\sigma}_{L^2(\I)}.
\]
\end{lemma}

\vspace{4mm}
\section{The energy functionals}\label{sec:ftionals}

We first establish a connection between the discrete energy functional $F^N$
and its ``semi-discrete" counterpart $\Fcont^N$.

\begin{lemma}\label{lem:reformEnergy}
For any $N\in\N$ and $B\in\BB^N_r$ we have
\[
   \Fcont^N[\QQ^N[B]] = F^N[B],
\]
where $F^N$ is the discrete energy functional given by \eqref{eq:F}
and $\Fcont^N$ is the integral energy functional \eqref{Fn:N}.
\end{lemma}

\begin{proof}
Let us denote $\b:=\QQ^N[B]$. The metabolic term in $\Fcont^N[\b]$ reads
\[
   \frac{\nu}{\gamma} \iintII \b(x,y)^\gamma \left[\ell^N(x,y) \right]^{\gamma+1} \d w^N(x, y) &=&
     \frac{\nu}{\gamma} \sum_{i=1}^N \sum_{j=1}^N W_{ij}^N B_{ij}^\gamma \left( \L_{ij}^N \right)^{\gamma+1} \int_{\I_i^N}\int_{\I_j^N} \d x\d y \\
       &=& \frac{\nu}{\gamma N^2} \sum_{i=1}^N \sum_{j=1}^N W_{ij}^N B_{ij}^\gamma \left( \L_{ij}^N \right)^{\gamma+1},
\]
where we used the fact that, by construction, $\b^N = \QQ^N[B] \equiv B_{ij}$ on the patch $\I_i^N\times \I_j^N$,
and similarly $\ell^N = \QQ^N[\L^N] \equiv \L^N_{ij}$, $w^N = \QQ^N[W^N] \equiv W^N_{ij}$ on $\I_i^N\times \I_j^N$.

For the kinetic term in $\Fcont^N[\b]$ we have
\[
    \iintII \b(x,y) \left(\dd p^N(x,y)\right)^2 \d w^N(x, y)
    = \sum_{i=1}^N \sum_{j=1}^N W_{ij}^N B_{ij} \int_{\I_i^N}\int_{\I_j^N} \left(\dd p^N(x,y)\right)^2 \d x\d y,
\]
where $p^N \in \LtwozeroN(\I)$ is the unique solution of the approximate Poisson equation \eqref{eq:Poissonsc} with permeability kernel $\b=\QQ^N[B]$.
For $(x,y) \in \I_i^N\times \I_j^N$ we have
\[
   \dd p^N(x,y) = p^N(y) - p^N(x)  = P_j - P_i,
\]
where we used \eqref{P_i}, $(P_i)_{i=1}^N$ being a solution of the Kirchhoff law \eqref{eq:KB} with the conductivities $(B_{i,j})_{i,j=1}^N$
and source/sink terms \eqref{ass:S}.
Consequently,
\[
   \Fcont^N[\b] =  \frac{1}{2N^2} \sum_{i=1}^N \sum_{j=1}^N \left( B_{ij} (P_i - P_j)^2 + \frac{\nu}{\gamma} B_{ij}^\gamma \left( \L_{ij}^N \right)^{\gamma+1} \right) W_{ij}^N = F^N[B].
\]
\end{proof}

In the sequel it proves advantageous to analyze the kinetic and metabolic parts of the energy functionals separately.
We start with the kinetic part,
\(  \label{def:Fkin}
      \Fkin[\b] :=  \frac12 \iintII \b \left(\dd p\right)^2 \d w(x,y),
\)
where $p \in \Ltwozero(\I)$ is the unique solution of the Poisson equation \eqref{eq:Poissonbw} with kernel $\b$, and
\[  
      \Fkin^N[\b] :=  \frac12 \iintII \b \left(\dd p^N\right)^2 \d w^N(x,y),
\]
with $p^N \in \LtwozeroN(\I)$ the unique solution of the approximate Poisson equation \eqref{eq:Poissonsc}.

\begin{lemma}\label{lem:ZZ}
For any $N\in\N$ and $\b\in L^2_r(\II)$ we have
\[
   \Fkin^N[\b] = \Fkin^N[\ZZ^N[\b]].
\]
\end{lemma}

\begin{proof}
Let $p^N \in \LtwozeroN(\I)$ be the unique weak solution of the approximate Poisson equation \eqref{eq:Poissonsc}
with kernel $\b$, i.e.,
\(   \label{eq:ZZZ}
   \frac12 \iintII \b(x,y) (\dd p^N(x,y)) (\dd \varphi^N(x,y)) \d w^N(x,y) = \int_\I \sigma(x) \varphi^N(x) \d x
\)
for any test function $\varphi^N\in L_{0,N}^\infty(\I)$.
By definition, $\varphi^N$, $p^N$ and $w^N$ are piecewise constant, so that we have
$\varphi^N = \ZZ^N[\varphi^N]$ and analogously for $p^N$ and $w^N$. Consequently,
\[
   (\dd p^N) (\dd \varphi^N) w^N = \ZZ^N\left[ (\dd p^N) (\dd \varphi^N) w^N \right].
\]
Inserting this in the left-hand side of \eqref{eq:ZZZ} and using the self-adjoint property of $\ZZ^N$ on $L^2(\II)$, we obtain
\[
   \frac12 \iintII  \ZZ^N[\b](x,y) (\dd p^N(x,y)) (\dd \varphi^N(x,y)) \d w^N(x,y) = \int_\I \sigma(x) \varphi^N(x) \d x.
\]
We conclude that $p^N$ is the unique weak solution of \eqref{eq:Poissonsc}
with kernel $\ZZ^N[\b]$. Then identity \eqref{eq:testPoisson-h} of Lemma \ref{lem:Poisson-bh}
gives
\[
   \Fkin^N[\b] &=& \frac12 \iintII \b(x,y) (\dd p^N(x,y))^2 \d w^N(x,y)  \\
   &=& \int_\I \sigma(x) p^N(x) \d x \\
   &=& \frac12 \iintII  \ZZ^N[\b](x,y) (\dd p^N(x,y))^2 \d w^N(x,y) \\
   &=& \Fkin^N[\ZZ^N[\b]].
\]
\end{proof}

Next we prove strong continuity of the kinetic energy functional $\Fkin$ on $L^2(\II)$.

\begin{lemma}\label{lem:strong}
Let $\sigma \in \Ltwozero(\I)$ and $r>0$.
Let the sequence $(\b^n)_{n\in\N} \subset L^2_r(\II)$
converge to $b\in L^2_r(\II)$ in the norm topology of $L^2(\II)$ as $n\to \infty$.
Then
\[
   \Fkin[b] = \lim_{n\to\infty} \Fkin[b^n].
\]
\end{lemma}

\begin{proof}
Let $(p^n)_{n\in\N}\subset \Ltwozero(\I)$ be a sequence of
weak solutions of the Poisson equation \eqref{eq:Poissonbw} with permeability kernel $\b^n\in L^2_r(\II)$,
provided by Lemma \ref{lem:Poisson-b}.
Identity \eqref{eq:testPoisson} gives
\(  \label{eq:testPoisson2}
   \frac12 \iintII \b^n |\dd p^n|^2 \d w(x,y) = \int_\I \sigma p^n \d x\d y.
\)
Due to the uniform bound \eqref{est:p},
there exists $p\in \Ltwozero(\I)$, a weak limit of (a subsequence of) $p^n$,
verifying \eqref{eq:Poissonbw} for all test functions $\varphi\in L_0^\infty(\I)$.
Since $\b^n\to \b$ strongly in $L^2(\II)$, we can pass to the limit in the weak formulation \eqref{eq:Poissonbw}
which establishes $p$ as its weak solution with permeability kernel $\b$
and test functions $\varphi\in L_0^\infty(\I)$.
Again, \eqref{eq:testPoisson} gives
\[
   \frac12 \iintII \b (\dd p)^2 \d x\d y = \intI \sigma p \d x.
\]
Combining with \eqref{eq:testPoisson2} we have
\[
   \lim_{n\to\infty} \frac12 \iintII b^n (\dd p^n)^2 \d x\d y =
   \lim_{n\to\infty}  \intI \sigma p^n \d x =  \intI \sigma p \d x =
   \frac12 \iintII b (\dd p)^2 \d x\d y,
\]
which is the strong continuity of the kinetic part of $\Fcont[\b^n]$.
\end{proof}

Our next result shows that the functional $\Fkin$ is convex on its domain.
Convexity of the kinetic part of the discrete energy functional \eqref{en:E} and its continuum counterpart
has been established in \cite[Proposition 3.2]{HMP:CMS:2022} and in \cite[Section 3]{HMPort:23},
based on a proof of positive semidefinitness of the Hessian.
Here we offer a rigorous proof using the variational formulation of the Poisson equation \eqref{eq:Poissonbw}.

\begin{lemma}\label{lem:convex}
The functional $\Fkin$ defined in \eqref{def:Fkin} 
is finite and convex on the convex set $L^2_r(\II)$.
\end{lemma}

\begin{proof}
Finiteness of $\Fkin$ on $L^2_r(\II)$ follows directly from \eqref{eq:testPoisson}. To prove convexity on $L^2_r(\II)$, we first prove convexity on $L_r^\infty(\II)$. According to Lemma~\ref{lem:Poisson}, for $b \in L_r^\infty(\II)$, the Poisson equation \eqref{eq:Poissonbw} is well posed in $L^2_0(\I)$ with test functions in $L^2_0(\I)$. Moreover, its unique solution $p=p[b]\in L^2_0(\I)$ is the unique minimizer of the functional
\[
\mathcal{J}[v;b]:=\frac14\iint_{\II} b(x,y)\,(\dd v(x,y))^2\,\d w(x,y)-\int_\I \sigma(x)\,v(x)\,\d x
\qquad \text{over } v\in L^2_0(\I).
\]
Indeed, \eqref{eq:Poissonbw} is the Euler--Lagrange equation of $\mathcal{J}[\cdot;b]$.
Testing \eqref{eq:Poissonbw} with $\varphi=p$ yields
$\frac12\iint_{\II} b(\dd p)^2\,\d w=\int_\I \sigma p\,\d x$, hence
\[
\mathcal{J}[p;b]=\frac14\iint_{\II} b(\dd p)^2\,\d w-\int_\I \sigma p\,\d x
= \frac12\,\Fkin[b]-\Fkin[b]=-\frac12\,\Fkin[b],
\]
and, consequently,
\begin{equation}\label{eq:dualFkin}
\Fkin[b]= -2\min_{v\in L^2_0(\I)} \mathcal{J}[v;b]
= \sup_{v\in L^2_0(\I)}\left\{2\int_\I \sigma v\,\d x-\frac12\iint_{\II} b(\dd v)^2\,\d w\right\}.
\end{equation}
For each fixed $v\in L^2_0(\I)$ the map
\[
b\ \mapsto\ 2\int_\I \sigma(x) v(x)\,\d x-\frac12\iintII b(x,y)(\dd v(x,y))^2\,\d w(x,y)
\]
is affine in $b$. Therefore \eqref{eq:dualFkin} expresses $\Fkin$ as a supremum of affine functions of $b$, which implies that $\Fkin$ is convex on $L_r^\infty(\II)$.
 \par
Now, let $b_0,b_1\in L^2_r(\II)$ and $\theta\in[0,1]$. For $M\ge r$ define the truncations
\[
b_i^M:=\min\{b_i,M\}\in L_r^\infty(\II),\qquad i\in\{0,1\}.
\]
Set $b_\theta:=(1-\theta)b_0+\theta b_1$ and $b_\theta^M:=(1-\theta)b_0^M+\theta b_1^M$. We have, for every $M\ge r$,
\[
\Fkin[b_\theta^M]\ \le\ (1-\theta)\Fkin[b_0^M]+\theta \Fkin[b_1^M].
\]
Moreover, $b_i^M\to b_i$ strongly in $L^2(\II)$ as $M\to\infty$, hence also $b_\theta^M\to b_\theta$ strongly in $L^2(\II)$.
By the strong continuity of $\Fkin$ on $L^2_r(\II)$ established in Lemma~\ref{lem:strong} we conclude
\begin{align*}
\Fkin[b_\theta] &=\lim_{M\to\infty}\Fkin[b_\theta^M] \\
&\le (1-\theta)\lim_{M\to\infty}\Fkin[b_0^M]+\theta\lim_{M\to\infty}\Fkin[b_1^M] \\
&=(1-\theta)\Fkin[b_0]+\theta\Fkin[b_1].
\end{align*}
This proves convexity of $\Fkin$ on $L^2_r(\II)$.
\end{proof}

We now observe that $\Fkin$ is strongly continuous on $L^2_r(\II)$
due to Lemma \ref{lem:strong},
and convex on the same set due to Lemma \ref{lem:convex}.
Then, as a consequence of Mazur's lemma \cite{Rockafellar}, it is weakly lower semicontinuous, i.e.,
\(  \label{eq:wlsc}
   \Fkin[\b] \leq \liminf_{n\to \infty} \Fkin[\b^n]
\)
whenever $(\b^n)_{n\in\N} \subset L^2_r(\II)$ converges weakly in $L^2(\II)$ to $\b$.

\vspace{4mm}
\section{Proof of Theorem \ref{thm:1}}\label{sec:thm:1}

We organize the proof into the following four lemmas, working separately with the kinetic and metabolic parts of the energy functional.
We start by proving the $\liminf$-part of the $\Gamma$-convergence, i.e., claim \eqref{eq:gamma:liminf} of Theorem \ref{thm:1}.

\begin{lemma}\label{lem:liminf_kin}
Let the sequence $\b^N\in L^2_r(\II)$ converge weakly in $L^2(\II)$ to $\b\in L^2_r(\II)$.
Then
\[  
   \Fkin[b] \leq  \liminf_{N\to\infty}  \Fkin^N[\b^N].
\]
\end{lemma}

\begin{proof}
We observe that the weak convergence $\b^N\wto\b$ in $L^2(\II)$
gives $\ZZ^N[\b^N] \wto\b$ in the same space.
Indeed, for any test function $\varphi\in L^2(\II)$ we have,
using the self-adjoint property of the orthogonal projection $\ZZ^N$ on $L^2(\I)$,
\[
   \iintII \ZZ^N[\b^N] \varphi \d x\d y &=& \iintII \b^N \ZZ^N[\varphi] \d x\d y \\
     &=&  \iintII \b^N \left( \ZZ^N[\varphi] - \varphi\right) \d x \d y + \iintII \b^N \varphi \d x \d x,
\]
and the first integral on the right-hand side vanishes in the limit $N\to\infty$
due to the uniform boundedness of $b^N$ in $L^2(\II)$
and the approximation property
\(  \label{eq:approxProp}
   \lim_{N\to\infty} \Norm{ \ZZ^N[\varphi] - \varphi}_{L^2(\II)} = 0.
\)
Therefore, the weak lower semicontinuity \eqref{eq:wlsc} of the functional $\Fkin$ gives
\(   \label{wlscG}
    \Fkin[\b]  \leq \liminf_{N\to\infty} \Fkin[\ZZ^N[\b^N] ]. 
\)
Identity \eqref{eq:testPoisson} of Lemma \ref{lem:Poisson-b} gives
\[
   \Fkin[\ZZ^N[\b^N] ] =  \frac12 \iintII \ZZ^N[\b^N]  \left(\dd p^N\right)^2 \d w(x,y) = \intI \sigma p^N \d x,
\]
where $p^N\in \Ltwozero(\I)$ is the unique solution of the Poisson equation \eqref{eq:Poissonbw} with permeability kernel $\ZZ^N[\b^N]$.
Analogously, using identity \eqref{eq:testPoisson-h} of Lemma \ref{lem:Poisson-bh},
\[
   \Fkin^N[\ZZ^N[\b^N] ] =  \frac12 \iintII \ZZ^N[\b^N] \left(\dd p^N_h\right)^2 \d w^N(x, y) = \intI \sigma p_h^N \d x,
\]
where $p^N_h\in \LtwozeroN(\I)$ is the unique solution of the semidiscrete Poisson equation \eqref{eq:Poissonsc} with permeability kernel $\ZZ^N[\b^N]$.
Subtracting the above two identities, we have
\(   \label{eq:FkinGkin}
   \Fkin[\ZZ^N[\b^N]] - \Fkin^N[\ZZ^N[\b^N] ] = \intI \sigma (p^N - p_h^N) \d x.
\)
Since $\ZZ^N[\b^N]$ is piecewise constant, identity \eqref{P_i} of Lemma \ref{lem:Poisson} gives $p_h^N = \ZZ^N[p^N]$.
Using again the self-adjointness of $\ZZ^N$ on $L^2(\I)$, we have
\[
    \intI \sigma (p^N - \ZZ^N[p^N]) \d x = \intI (\sigma - \ZZ^N[\sigma]) p^N \d x,
\]
and the right-hand side vanishes in the limit $N\to\infty$ due to the uniform bondedness of $p^N$ in $L^2(\I)$
and the approximation property \eqref{eq:approxProp} for $\sigma$.
Recalling \eqref{wlscG}, we conclude that
\[
   \Fkin[\b]  &\leq& \liminf_{N\to\infty} \Fkin[\ZZ^N[\b^N] ]\\
    &=& \liminf_{N\to\infty} \left[ \Fkin^N[\ZZ^N[\b^N] ] + \intI \sigma (p^N - p_h^N) \d x \right]  \\
     &=& \liminf_{N\to\infty}  \Fkin^N[\b^N]  + \lim_{N\to\infty} \intI (\sigma - \ZZ^N[\sigma]) p^N  \d x \\
     &=& \liminf_{N\to\infty}  \Fkin^N[\b^N],
\]
where we used the identity $\Fkin^N[\ZZ^N[\b^N]] = \Fkin^N[\b^N]$ provided by Lemma \ref{lem:ZZ}.
\end{proof}

\begin{lemma}\label{lem:liminf_met}
Let $\gamma>1$ and let the sequence $\b^N\in L^\gamma(\II)$ converge weakly in $L^\gamma(\II)$ to $\b\in L^\gamma(\II)$.
Then
\[   
   \Fmet[b] \leq  \liminf_{N\to\infty}  \Fmet^N[\b^N].
\]

\end{lemma}

\begin{proof}
Let us write the metabolic part of $\Fcont^N$ as
\[
   \Fmet^N[\b^N] &=& \frac{\nu}{2\gamma} \iintII \left( \b^N(x,y) \right)^\gamma \left( \ell^N(x,y) \right)^{\gamma+1} \d w^N(x, y)
   = \frac{\nu}{2\gamma} \Norm{a^N}^\gamma_{L^\gamma(\II)},
\]
where we denoted
\[
   a^N := \b^N \left(\ell^N\right)^\frac{\gamma+1}{\gamma} w^N.
\]
Since $\Norm{\ell^N}_{L^\infty(\II)} \leq 1$ for all $N\in\N$ and $w^N$ only has values in $\{0,1\}$,
we obviously have $a^N$ uniformly bounded in $L^\gamma(\II)$.
Therefore, up to an eventual extraction of a subsequence, it converges weakly to some $a\in L^\gamma(\II)$.
The weak lower semicontinuity of the $L^\gamma$ norm readily gives
\(  \label{eq:norma}
   \Norm{a}^\gamma_{L^\gamma(\II)} \leq \liminf_{N\to\infty} \Norm{a^N}^\gamma_{L^\gamma(\II)} = \frac{2 \gamma}{\nu} \liminf_{N\to\infty} \Fmet^N[\b^N].
\)
To identify the limit, we write for a test function $\varphi\in L^\infty(\II)$,
\(
   \label{eq:identify}
   \iintII a^N\varphi \d x\d y &=& \iintII b^N \left(\ell^N\right)^\frac{\gamma+1}{\gamma} \varphi \,\d w^N(x, y) \\
   &=& \iintII b^N \left( \left(\ell^N\right)^\frac{\gamma+1}{\gamma} - \ell^\frac{\gamma+1}{\gamma} \right) \varphi \d w^N(x, y) + 
     \iintII b^N \ell^\frac{\gamma+1}{\gamma} \varphi \, \d w^N(x, y).
     \nonumber
\)
The first term on the right-hand side vanishes in the limit $N\to\infty$ due to the uniform boundedness of $b^N w^N$
in $L^\gamma(\II)$ and the strong convergence of $\ell^N$ to $\ell$ in $L^q(\II)$ with any $q<\infty$, imposed by assumption \eqref{ass:L2}.
The second term converges to $\iintII b \,\ell^\frac{\gamma+1}{\gamma} \varphi \,\d w(x, y)$
due to the weak convergence $\b^N \wto b$ in $L^\gamma(\II)$
and strong convergence of $w^N$ to $w$ in $L^q(\II)$ with any $q<\infty$, given by assumption \eqref{ass:A2}.
We conclude that $a = \b \ell^\frac{\gamma+1}{\gamma} w$ and with \eqref{eq:norma} we arrive at
\[
   \Fmet[b] = \frac{\nu}{2\gamma} \Norm{a}^\gamma_{L^\gamma(\II)} \leq \liminf_{N\to\infty} \Fmet^N[\b^N].
\]
\end{proof}

Combining the claims of Lemmas \ref{lem:liminf_kin} and \ref{lem:liminf_met}, we conclude
that for any sequence $\b^N\in L_r^\omega(\II)$, $\omega:=\max\{2,\gamma\}$,
converging weakly in $L^\omega(\II)$ to $\b\in L^\omega(\II)$, we have
\[
   \Fcont[\b] &=& \Fkin[\b] + \Fmet[\b] \\
      &\leq& \liminf_{N\to\infty}  \Fmet^N[\b^N] + \liminf_{N\to\infty}  \Fkin^N[\b^N] \\
      &\leq& \liminf_{N\to\infty} \left( \Fmet^N[\b^N] + \Fkin^N[\b^N] \right) \\
      &=& \liminf_{N\to\infty}  \Fcont^N[\b^N].
\]
This proves claim \eqref{eq:gamma:liminf} of Theorem \ref{thm:1}.

\begin{lemma}\label{lem:limsup_kin}
Let $\b\in L^2(\II)$.
Then the sequence $\b^N := \ZZ^N[\b]$ converges in the norm topology of $L^2(\II)$ to $\b$ as $N\to\infty$, and
\[  
   \Fkin[b] =  \lim_{N\to\infty}  \Fkin^N[\b^N].
\]
\end{lemma}

\begin{proof}
The convergence of $\b^N := \ZZ^N[\b]$ to $\b$ strongly in $L^2(\II)$ is a standard result of the approximation theory.
The strong continuity of the functional $\Fkin$ established in Lemma \ref{lem:strong}, gives
\[
    \Fkin[\b]  = \lim_{N\to\infty} \Fkin[\b^N]. 
\]
Moreover, similarly as in \eqref{eq:FkinGkin},
\[
   \Fkin[\b^N] - \Fkin^N[\b^N] = \intI \sigma (p^N - p_h^N) \d x,
\]
where $p^N\in \Ltwozero(\I)$ is the unique solution of the Poisson equation \eqref{eq:Poissonbw},
and $p^N_h\in \LtwozeroN(\I)$ is the unique solution of the semidiscrete Poisson equation \eqref{eq:Poissonsc}, both with permeability kernel $\b^N$.
Since $\b^N = \ZZ^N[\b]$ is piecewise constant, by identity \eqref{P_i} of Lemma \ref{lem:Poisson} we have $p_h^N = \ZZ^N[p^N]$.
Consequently,
\[
    \intI \sigma (p^N - \ZZ^N[p^N]) \d x = \intI (\sigma - \ZZ^N[\sigma]) p^N \d x.
\]
The right-hand side vanishes in the limit $N\to\infty$ due to the uniform bondedness of $p^N$ in $L^2(\I)$
provided by \eqref{est:p}.
We thus have
\[
   \Fkin[\b]  &=& \lim_{N\to\infty} \left[ \Fkin^N[\b^N] + \intI (\sigma - \ZZ^N[\sigma]) p^N \d x \right] \\
   &=& \lim_{N\to\infty} \Fkin^N[\b^N].
\]
\end{proof}

\begin{lemma}\label{lem:limsup_met}
Let $\b\in L^\gamma(\II)$.
Then the sequence $\b^N := \ZZ^N[\b]$ converges in the norm topology of $L^\gamma(\II)$ to $\b$ as $N\to\infty$, and
\[  
   \Fmet[b] =  \lim_{N\to\infty}  \Fmet^N[\b^N].
\]
\end{lemma}

\begin{proof}
We have
\[
   \lim_{N\to\infty} \Fmet^N[\b^N] &=& \lim_{N\to\infty}  \iintII \left( \b^N(x,y) \right)^\gamma \left( \ell^N(x,y) \right)^{\gamma+1} \d w^N(x, y)  \\
      &=& \iintII \left( \b(x,y) \right)^\gamma \left( \ell(x,y) \right)^{\gamma+1} \d w(x, y)  \\
      &=& \Fmet[\b]
 \]
due to the strong convergence $\b^N\to \b$ in $L^\gamma(\II)$ and $\ell^N \to \ell$, $w^N \to w$ in $L^q(\II)$ with any $q<\infty$
and the uniform boundedness $\Norm{\ell^N w^N}_{L^\infty(\II)} \leq 1$.
\end{proof}

Combining the claims of Lemmas \ref{lem:limsup_kin} and \ref{lem:limsup_met}, we conclude
the for any $\b\in L^\omega(\II)$, $\omega:=\max\{2,\gamma\}$, the sequence $\b^N := \ZZ^N[\b]$
converges to $b$ strongly in $L^\omega(\II)$ and
\[
   \Fcont[\b] = \lim_{N\to\infty} \Fcont^N[\b^N].
\]
This proves claim \eqref{eq:gamma:limsup} of Theorem \ref{thm:1}.

\vspace{4mm}
\section{Proof of Theorem \ref{thm:2}}\label{sec:thm:2}

\begin{proof}
Let $B^N \in \BB^N_r$, $N\in\N$, be a sequence of global minimizers of 
the discrete energy functionals $F^N$ given by \eqref{eq:F}--\eqref{eq:KB}.
We then have
\begin{equation}\label{eq:est_r}
   F^N[B^N] \leq F^N[r] = \frac{1}{2N^2} \sum_{i=1}^N \sum_{j=1}^N 
\left( (P_j-P_i)^2 r + \frac{\nu}{\gamma} r^{\gamma}  
(\L^N_{ij})^{\gamma+1} \right) W_{ij}^N,
\end{equation}
where $(P_i)_{i=1}^N$ is a solution of the Kirchhoff law \eqref{eq:KB} 
with conductivities $r>0$.
Since $\L^N_{ij} \leq 1$ and $W^N_{ij}\in\{0,1\}$, we have
\[
   \frac{1}{2N^2} \sum_{i=1}^N \sum_{j=1}^N \frac{\nu}{\gamma} 
r^{\gamma}  (\L^N_{ij})^{\gamma+1}  W_{ij}^N \leq \frac{\nu}{2\gamma} 
r^{\gamma}.
\]
Moreover, estimate \eqref{est:Kirchhoff} of Lemma \ref{lem:Kirchhoff} gives
\[
   \frac{r}{2N^2} \sum_{i=1}^N \sum_{j=1}^N  (P_j-P_i)^2 W_{ij}^N  \leq 
\frac{4}{r \lambda^2} \intI \sigma(x)^2 \d x.
\]
Consequently, the right-hand side in \eqref{eq:est_r} is uniformly bounded.
Denoting $\b^N := \QQ^N[B^N]$, Lemma \ref{lem:reformEnergy} gives
$\Fcont^N[\b^N] = F^N[B^N]$, and we have
\[
   \sup_{N\in\N} \Fkin^N[\b^N] <\infty \qquad\mbox{and}\qquad  
\sup_{N\in\N} \Fmet^N[\b^N] <\infty.
\]
We define the sequence $a^N := \b^N \left(\ell^N\right)^\frac{\gamma+1}{\gamma} w^N$. 
Note that the metabolic energy can be written as
\[
   \Fmet^N[\b^N] = \frac{\nu}{2\gamma}\iintII \left( \b^N \right)^\gamma \left( \ell^N \right)^{\gamma+1} \d w^N(x,y) = \frac{\nu}{2\gamma} \Norm{a^N}^\gamma_{L^\gamma(\II)}.
\]
The uniform bound on $\Fmet^N[\b^N]$ implies that the sequence $a^N$ is uniformly bounded in $L^\gamma(\II)$.
Since $\gamma > 1$, the space $L^\gamma(\II)$ is reflexive, and there exists a subsequence (not relabeled) such that $a^N \wto a$ weakly in $L^\gamma(\II)$.
Using the weak lower semicontinuity of the $L^\gamma$-norm, we have
\[
   \Norm{a}^\gamma_{L^\gamma(\II)} \leq \liminf_{N\to\infty} 
\Norm{a^N}^\gamma_{L^\gamma(\II)} = \frac{2 \gamma}{\nu}\liminf_{N\to\infty} \Fmet^N[\b^N].
\]
We identify the limit $a\in L^\gamma(\II)$. Since $\b^N=\QQ^N[B^N]$ and $B^N_{ij}=0$ whenever $W^N_{ij}=0$, we have $\b^N=\b^N w^N$, and thus
$a^N=\b^N \left(\ell^N\right)^\frac{\gamma+1}{\gamma}$.
Consequently,
\[
   \b^N = a^N \left(\ell^N\right)^{-\frac{\gamma+1}{\gamma}} .
\]
Since $\gamma>2$, by H\"older's inequality we have
\[
   \Norm{\b^N}_{L^2(\II)}
   \le \Norm{a^N}_{L^\gamma(\II)} \Norm{\left(\ell^N\right)^{-\frac{\gamma+1}{\gamma}}}_{L^{\frac{2\gamma}{\gamma-2}}(\II)}
   = \Norm{a^N}_{L^\gamma(\II)} \Norm{\left(\ell^N\right)^{-1}}_{L^{\frac{2(\gamma+1)}{\gamma-2}}(\II)}^{\frac{\gamma+1}{\gamma}},
\]
and the right-hand side is uniformly bounded by assumption~\eqref{ass:L3star}. Therefore, there exists $b\in L^2(\II)$ and a subsequence (not relabeled) such that
$\b^N \wto b$ weakly in $L^2(\II)$. \\
Now, for any $\varphi\in L^\infty(\II)$ we have
\[
   \iintII a^N \varphi \d x\d y
   = \iintII \b^N \left(\ell^N\right)^\frac{\gamma+1}{\gamma} w^N \varphi \d x\d y .
\]
We write
\begin{align*}
\iintII \b^N \left(\ell^N\right)^\frac{\gamma+1}{\gamma} w^N \varphi \d x\d y
&= \iintII \b^N \ell^\frac{\gamma+1}{\gamma} w \varphi \d x\d y \\
&\quad + \iintII \b^N \Big( \left(\ell^N\right)^\frac{\gamma+1}{\gamma}-\ell^\frac{\gamma+1}{\gamma} \Big) w^N \varphi \d x\d y \\
  &\quad + \iintII \b^N \ell^\frac{\gamma+1}{\gamma} (w^N-w)\varphi \d x\d y .
\end{align*}
Since $\ell^N\to \ell$ in $L^1(\II)$ by assumption \eqref{ass:L2}, and $\Norm{\ell^N}_{L^\infty(\II)}\le 1$,
the map $t\mapsto t^{\frac{\gamma+1}{\gamma}}$ is Lipschitz continuous on~$[0,1]$, and thus
$\left(\ell^N\right)^\frac{\gamma+1}{\gamma}\to \ell^\frac{\gamma+1}{\gamma}$ in $L^1(\II)$.
Moreover, $w^N\to w$ in $L^1(\II)$ and $\Norm{w^N}_{L^\infty(\II)}\le 1$.
Therefore,
\[
\Big\|\Big( \left(\ell^N\right)^\frac{\gamma+1}{\gamma}-\ell^\frac{\gamma+1}{\gamma} \Big) w^N\Big\|_{L^2(\II)}\to 0,
\qquad
\|w^N-w\|_{L^2(\II)}\to 0,
\]
as $N\to\infty$.
Using the uniform boundedness of $\b^N$ in $L^2(\II)$,
we conclude
\[
   \lim_{N\to\infty} \iintII \b^N \left(\ell^N\right)^\frac{\gamma+1}{\gamma} w^N \varphi \d x\d y &=&
   \lim_{N\to\infty} \b^N \ell^\frac{\gamma+1}{\gamma} w \varphi \d x\d y \\
    &=& \iintII b \ell^\frac{\gamma+1}{\gamma} w \varphi \d x\d y,
\]
where we used the weak convergence $\b^N \wto b$ in $L^2(\II)$.
We thus have $a = b \ell^\frac{\gamma+1}{\gamma} w$,
and, consequently,
\begin{equation}\label{eq:Fmet}
\begin{split}
\Fmet[b]
& = \frac{\nu}{2 \gamma} \iint_{\II} b(x,y)^\gamma \,\ell(x,y)^{\gamma+1}\,\d w(x,y) \\
& = \frac{\nu}{2 \gamma} \|a\|_{L^\gamma(\II)}^\gamma \\
& \le \liminf_{N\to\infty} \Fmet^N[\b^N].
\end{split}
\end{equation}
Moreover, Lemma~\ref{lem:liminf_kin} gives
\[
\Fkin[b] \le \liminf_{N\to\infty} \Fkin^N[\b^N],
\]
and combining this with \eqref{eq:Fmet} we obtain
\[
\Fcont[b] \le \liminf_{N\to\infty} \Fcont^N[\b^N].
\]

We claim that $b$ is a global minimizer of the continuum energy functional \eqref{eq:Fcont}--\eqref{eq:Poissonbw} in $L^\gamma_r(\II)$.
For contradiction, let us assume that there exists $\barb\in L^\gamma_r(\II)$ such that
\begin{equation}\label{eq:contra}
   \Fcont[\barb] < \Fcont[b].
\end{equation}
Then Theorem \ref{thm:1} asserts that
the sequence $\barb^N := \ZZ^N[\barb]$ converges to $\barb$ in the norm topology of $L^\gamma(\II)$,
and \eqref{eq:gamma:limsup} gives
\[
   \Fcont[\barb] = \lim_{N\to\infty} \Fcont^N[\barb^N].
\]
However, since each $B^N$ is, per construction, a global minimizer of $F^N$,
we have
\[
   \Fcont^N[\b^N] = F^N[B^N]  \leq F^N[\tilde B^N] = \Fcont^N[\barb^N],
\]
where $\tilde B^N \in \BB^N_r$ is the matrix of the values of $\barb^N$, i.e., $\barb^N = \QQ^N[\tilde B^N]$.
Therefore,
\[
     \Fcont[\barb] = \lim_{N\to\infty} \Fcont^N[\barb^N] \geq \liminf_{N\to\infty} \Fcont^N[\b^N]  \geq \Fcont[b],
\]
which is a contradiction to \eqref{eq:contra}.
\end{proof}

\vspace{4mm}
\section*{Acknowledgments}
This work was partially supported by the Austrian Science Fund (FWF), project number 10.55776/F65.



\end{document}